\documentclass[12pt]{amsart}
\usepackage{amsmath}
\newcommand{\C}{\mathbb C}

\newcommand{\Fn}{{\mathcal F}_n}
\newtheorem{thm}{Theorem}
\newtheorem{lem}{Lemma}

\begin{document}
\title{A Proof of Smale's Mean Value Conjecture}
\author{Gerald Schmieder}
\date{}
\begin{abstract} A proof of Smale's mean value conjecture from 1981 is
given.
\end{abstract}
\maketitle Connected with his investigations on the complexity of
determining polynomial roots by Newton's method, Steve Smale
\cite{smale} considered difference quotients
$$D(\zeta,z)=\frac{p(\zeta)-p(z)}{\zeta-z},$$ where $p$ is a
non-constant polynomial, $p'(\zeta)=0$ and $z\neq \zeta$ an
arbitrary complex number. He asked for some universal (i.e. valid
for all such polynomials and all $z\neq \zeta$) constant $K$ such
that $|D(\zeta,z)|\le K|p'(z)|$ for at least one derivative zero
$\zeta$. He proved in \cite{smale}, using results on univalent
functions, that this is true for $K=4$ and conjectured $K=1$ to be
best possible.\par Obviously one may without loss of generality
assume that $z=0$ and $p(0)=0$. Then the question is to estimate
the number $\displaystyle \min\big\{\big|\frac{p(\zeta)}{\zeta
p'(0)}\big|\,:\,p'(\zeta)=0\big\}$. Note that the conjecture
trivially holds for polynomials of degree one.
\smallskip

The conjectured bound $1$ can be sharpend a little bit if we
consider only polynomials of a fixed degree. Here we will prove
the following:\newline {\em Let $p \in \C\, [z]$ be a polynomial
of degree $n>1$ with $p(0)=0$ and $p'(0)\neq 0$. Then
$$\min\left\{\Big|\frac{p(\zeta)}{\zeta
p'(0)}\Big|\,:\,p'(\zeta)=0\right\}\le\frac{n-1}{n}.$$ Equality
only occurs for $p(z)=a_1 z + a_nz^n$ with arbitrary
$a_1,a_n\in\C\setminus\{0\}$.}

\smallskip

Let $n>1$ be fixed and define $\Fn$ as the class of $n$th degree monic complex
polynomials $p$ with $p(0)=0$, $p'(0)\neq 0$ and $p(\zeta)\neq 0$ for all
derivative zeros $\zeta$ of $p$. Obviously it suffices to consider polynomials
$p\in\Fn$ in order to give a proof of Smale's conjecture. For such $p$ we
define $$\rho(p,\zeta):=\Big|\frac{p(\zeta)}{\zeta p'(0)}\Big|$$
and the {\em associated number} as
$$\rho(p):= \min\left\{\rho(p,\zeta)\,:\,p'(\zeta)=0\right\}.$$
The zero $\zeta_0$ of $p'$ is {\em essential} if
$$\rho(p)= \Big|\frac{p(\zeta_0)}{\zeta_0 p'(0)}\Big|.$$
Note that a polynomial may have more than one essential derivative zero.
We call $p\in\Fn$ {\em simple} if $p''(\zeta)\neq 0$ for all
essential derivative zeros $\zeta$ of $p$.

A polynomial $p_0\in \Fn$ is {\em maximal} if $\rho(p)\le
\rho(p_0)$ for all $p\in\Fn$. Below we will determine the maximal
polynomials in $\Fn$. In the following we will prove:
\begin{thm}\label{thm1} For each $p\in\Fn$ there exists some
$q\in\Fn$ which zeros $w_2,\dots,w_n\neq 0$ have the same modulus
and it holds $\rho(q)\ge \rho(p)$.\end{thm} In order to prove
theorem \ref{thm1} we may without loss of generality assume that
$|z_j|\le 1$ holds for the zeros $z_2,\dots,z_n\neq 0$ of $p$, and
equality is taken for at least one of them. Otherwise we consider
the polynomial $s^np(z/s)$ with $s=\max\{|z_2|,\dots,|z_n|\}$. The
associated number of this polynomial is the same than those of
$p$. Moreover we may assume that $|z_n|<1$.

\section{The basic idea\label{basic}}
If $p\in\Fn$ is a polynomial with the zeros $z_2,\dots,z_n$
besides $0$ and the derivative zero $\zeta$ with $p(\zeta)\neq 0$,
then
\begin{eqnarray}\label{form1}
\frac{p'}{p}(\zeta)=0=\frac{1}{\zeta}+\sum_{j=2}^n\frac{1}{\zeta-
z_j}.\end{eqnarray}

As explained we may provide that $|z_j|\le 1$ for $j=2,\dots,n$
and $|z_n|<1$. We let $z_2,\dots,z_{n-1}$ be fixed and vary $z_n$,
i.e., we consider the polynomials
\begin{eqnarray}\label{Q}
Q(z,u)=(z-u)\,z\prod
_{j=2}^{n-1}(z-z_j)=(z-u)\,q(z).\end{eqnarray}
\smallskip

We assume for the moment that $\zeta$ is a zero of $p'$, but not a
zero of $p''$. The implicit function theorem (cf. \cite{nar})
shows the existence of a holomorphic function $\zeta(u)$ with
$\zeta(z_n)=\zeta$ and $\displaystyle \frac{\partial Q}{\partial
z}(\zeta(u),u)\equiv 0$, defined in a neighborhood of $z_n$. If we
move $u$ along a path $\gamma$ in $\C$ starting in $\gamma(0)=z_n$
then we have an unrestricted analytic continuation of
$\zeta(\gamma(t))$ if $\frac{\partial^2Q}{\partial
z^2}(\zeta(\gamma(t)),\gamma(t))\neq 0$ for all $t$. If the path
would meet these exceptional points, we would have at least a
continuation of $\zeta(\gamma(t))$ which is at least continuous in
such points. Note that the values of $\zeta(\gamma(t))$ with
respect to this continuation move on the Riemann surface $R$,
which is defined by the equation $Q'(z,u)=0$ (derivative with
respect to $z$). We will discuss this surface in section 2.
\smallskip

It comes out ($Q'$ denotes the derivative of $Q$ with respect to
$z$)
\begin{eqnarray}\label{3a}
\frac{Q(\zeta(u),u)}{\zeta(u)Q'(0,u)}=\Big(\frac{\zeta(u)}{u}-1\Big)
\prod_{j=2}^{n-1}\Big(\frac{\zeta(u)}{z_j}-1\Big)\\\nonumber=
\Big(\frac{\zeta(u)}{u}-1\Big)q(\zeta(u))\prod_{j=2}^{n-1}z_j^{-1}=:g(u,\zeta(u)).
\end{eqnarray}
Note that $\ln \rho(Q(.,u),\zeta(u))=\Re \log g(u,\zeta(u))$.
\par Let $p\in\Fn$ and $\zeta$ be a (not necessarily essential) derivative zero of $p$.
As above let $0,z_2,\dots,z_n\in\overline{E}$ be the zeros of $p$
and $|z_n|<1$. If $\gamma:[0,1]\to \C$ is a path with
$\gamma(0)=z_n, \gamma(1)=u$ we see $$\frac{d}{dt}\ln
\rho(Q(.,\gamma(t)),\zeta(\gamma(t)))=\frac{d}{dt}\Re \log
g(\gamma(t),\zeta(\gamma(t)))=\Re
\frac{\frac{d}{dt}g(\gamma(t),\zeta(\gamma(t)))}
{g(\gamma(t),\zeta(\gamma(t)))}.$$ Note that $\zeta(\gamma(t))$
depends on the path $\gamma$. So we have $$\ln
\rho(Q(.,u),\zeta(u))-\ln \rho(p,\zeta)
\hspace{5cm}$$$$=\int_0^{\,1} \frac{d}{dt}\ln
f(\gamma(t),\zeta(\gamma(t)))\,dt=\Re\int_0^{\,1} \frac{d}{dt}\log
g(\gamma(t),\zeta(\gamma(t)))\,dt.$$ The integrand can be
calculated as
$$\frac{\frac{d}{dt}g(\gamma(t),\zeta(\gamma(t)))}{g(\gamma(t),\zeta(\gamma(t)))}
=\gamma'(t)\left(\frac{\zeta'(\gamma(t))\gamma(t)-\zeta(\gamma(t))}
{\big(\zeta(\gamma(t)-\gamma(t))\big)\gamma(t)}+\zeta'(\gamma(t))\frac{q'}{q}(\zeta(\gamma(t)))\right)\,,$$
and this leads to $$\ln \rho(Q(.,u),\zeta(u))-\ln
\rho(p,\zeta)=\Re\int_\gamma \frac{\zeta'(v)v-\zeta(v)}
{\big(\zeta(v)-v\big)v}+\zeta'(v)\frac{q'}{q}(\zeta(v))\,dv.$$ The
right hand side can be written as $$\Re \int_\gamma
\frac{-1}{\zeta(v)-v}\cdot\frac{\zeta(v)}{v}+\zeta'(v)
\Big(\frac{1}{\zeta(v)-v}+\frac{q'}{q}(\zeta(v))\Big)\,dv.$$ From
(\ref{form1}) we obtain
$$0=\frac{Q'(\zeta(v),v)}{Q(\zeta(v),v)}=\frac{1}{\zeta(v)}+\frac{1}{\zeta(v)-v}+\frac{q'}{q}(\zeta(v)).$$
It comes out
\begin{eqnarray}\nonumber
\ln \rho(Q(.,u),\zeta(u))=\ln \rho(p,\zeta)-\Re \int_\gamma
\frac{1}{
\zeta(v)-v}\cdot\frac{\zeta(v)}{v}+\frac{\zeta'(v)}{\zeta(v)}
\,dv\,,
\end{eqnarray}
and therefore
\begin{eqnarray}\label{int}
\rho(Q(.,u)),\zeta(u))=\hspace{75mm}\\\nonumber\rho(p,\zeta)\cdot\Big|\exp\left(-\int_\gamma
\frac{1}{
\zeta(v)-v}\cdot\frac{\zeta(v)}{v}+\frac{\zeta'(v)}{\zeta(v)}\,dv\right)\Big|
\\\nonumber=\rho(p,\zeta)\cdot\big|\frac{\zeta}{\zeta(u)}\big|\cdot\Big|\exp\left(-\int_\gamma \frac{1}{
\zeta(v)-v}\cdot\frac{\zeta(v)}{v}\,dv\right)\Big| .\end{eqnarray}

\section{The Riemann surface $R$}
The Riemann surface $R$ of the derivative zeros of $Q$ is given by the equation
\begin{eqnarray}\label{gl} Q'(w)=q(w)+(w-u)q'(w)=0.\end{eqnarray}
This (actually compact) manifold
$R$ consists of the points $w$ (which are the derivative zeros of
$Q(.,u)$, and the equation gives local uniformizations of $R$, if
the derivative of $u=\varphi(w):=w+\frac{q}{q'}(w)$ with respect
to $w$ does not vanish (note that these branch points are also
described by $\frac{\partial^2Q}{\partial z^2}(w,u)= 0$). So the
points $w$ where $2q'(w)^2=q(w)q''(w)$ are branch points of the
surface. this branch points play in fact no special role on the
Riemann surface, their appearance depend on the special local
coordinates, which are given by the defining equation (example:
the surface of the square root is defined by $w^2=u$ with $0$ as a
branch point; if we add this point, it is conformally equivalent
to the plane resp. $\overline{\C}$). They can actually added as
''normal'' points to the surface and have simply connected
neighborhoods on which local coordinates can be found.\par $R$, as
a compact surface, may be regarded as a $(n-1)$-sheeted covering
of $\overline{\C}$, and $\varphi$ gives a canonically projection $R\to
\overline{\C}$.\par We define
\begin{eqnarray}\label{deff}
f(u,\zeta(u)):=\frac{\zeta}{\zeta(u)}\exp\left(-\int_{\gamma_u}
\frac{1}{
\zeta(v)-v}\cdot\frac{\zeta(v)}{v}\,dv\right),\end{eqnarray} where
$\gamma_u:[0,1]\to \C$ with $\gamma_u,(0)=z_n, \gamma_u(1)=u$ and
$\zeta(\gamma_u(0))=\zeta_0$ (some fixed derivative zero of $p$),
$\zeta(\gamma_u(1))=\zeta(u)$. By (\ref{int}) we have
\begin{eqnarray}\label{Qf}
\rho(Q(.,u)),\zeta(u))=\rho(p,\zeta_0)\cdot
|f(u,\zeta(u))|.\end{eqnarray} $f$ is. up to isolated
singularities, a holomorphic function on $R$, because it has this
property in the local coordinate $u\in\C$ (the case $u=\infty$ we
discuss separately). The holomorphy is not obviously clear in the
following cases.
\begin{enumerate}
\item [(i)] $Q'(0,u_0)=\frac{\partial Q}{\partial z}(0,u_0)=0$, or
\item [(ii)] $Q(w_1,u_1)=0$ (this includes the case $u=\zeta(u)$), or
\item [(iii)] $2q'(w_2)^2=q(w_2)\cdot q''(w_2)$ (branch points)
\end{enumerate}
(in case of (i) or (ii) the polynomial $Q$ does not belong to the
class $\Fn$). We discuss this three cases.\par\vspace{3mm} Case
(i): The polynomial $p$ has only simple zeros. So $Q'(0,u_0)=0$ is
only possible if $u=0$. A direct calculation gives that
$\displaystyle \rho(Q(.,0))=\rho(Q(.,0),0)=\frac{1}{2}$. Thus $f$
is has a removable singularity in $u=0$ if $\zeta(u)=0$. If
$\varphi(w_0)=0$, but $w_0\neq 0$ (and thus is not essential for
$Q(.,0)$) then we see that $f$ has a pole in $w_0$, because
$\rho(Q(.,w),w)\to \infty $ if $w\to w_0$.\par\vspace{3mm} Case
(ii): The assumption implies that $Q(.,u_1)$ has a multiple zero
in the point $w_1$. This is only possible if $u_1$ is one of the
zeros $z_2,\dots,z_{n-1}$ of $p$ ($u=0$ has already been
discussed) and $u_1=w_1$. By the definition we see that
$\rho(Q(.,u_1),w_1)=0$ and $\rho(Q(.,u_1),w)>0$ if
$\varphi(w)=u_1$ and $w\neq w_1$. So these singularities of $f$
are removable. Moreover we have $\rho(Q(.,u_1),w_1)=0$ in this
case.\par\vspace{3mm} Case (iii): If $2q'(w_2)^2=q(w_2)q''(w_2)$,
then $w_2\notin \{0,z_2,\dots,z_{n-1}\}$, because $q$ has only
simple zeros in these points. (\ref{deff}) shows that $f$ is
bounded in a neighborhood of the branch point $w_2$ on $R$. Again
we conclude that $f$ has a removable singularity in this
case.\medskip

\par We summarize:\begin{lem} The
function $f$ as defined in (\ref{deff}) is meromorphic on the
Riemann surface $R':=\{w\in R\,:\,\varphi(w)\in\C\}$. It has poles
exactly in the points $w\in R'$ with $\varphi(w)=0$ and $w\neq 0$.
The zeros of $f$ are the points $w\in R'$ with
$w=\varphi(w)\in\{z_2,\dots,z_{n-1}\}$.\end{lem} We can give an
alternative representation of $f$. It holds
$\rho(Q(.,u),\zeta(u))=\big|\frac{Q(\zeta(u),u)}{\zeta(u)Q'(0,u)}\big|$.
From (\ref{Qf}) we obtain that $f(u,\zeta(u))$ equals
$\frac{\zeta_0
p'(0)}{p(\zeta_0)}\frac{Q(\zeta(u),u)}{\zeta(u)Q'(0,u)}$, up to a
possible factor of modulus one. For $u=z_n$ we see that this
factor is one. By (\ref{gl}) we receive the representation:
\begin{eqnarray}\label{fneu}
f(u,\zeta(u))= \frac{z_n\, \zeta_0 \,q'(\zeta_0)}{q(\zeta_0)^2}
\cdot\frac{q(\zeta(u))^2}{u\,\zeta(u)\,q'(\zeta(u))}.
\end{eqnarray}
Finally we investigate the structure of $R$ close to $u=\infty$.
The point infinity is no branch point of $R$, because the function
$1/\varphi(1/w)$ has in $w=0$ the expansion
$w(\frac{n-1}{n}+a_1w+\dots)$.
\par For $u\in \overline{E}$ all zeros of $Q(.,u)$ are contained in
$\overline{E}$. By the Gau\ss{}-Lucas theorem we know that the
zeros of the derivative $Q'(z,u)=\frac{\partial Q}{\partial
z}(z,u)$ lie in the convex hull $C$ of the zeros. They are inner
points of $C$ with the only exception of multiple zeros of $Q$.
None of these derivative zeros in our case is of bigger order than
$1$. So the same argument gives that the zeros of the second order
derivative $Q''(z,u)=\frac{\partial^2Q}{\partial z^2}(z,u)$ are
points the open unit disk $E$. So the same is true for the branch
points of $R$. To be more precise, all branch points $w$ of $R$
fulfill $|\varphi(w)|<1$.\par The subset $D_1$ of $R$ with
$\varphi(D_1)=\overline{E}$ therefore contains all branch
points.\par As a consequence, the complement $R\setminus D_1$
(including $\infty$) consists of $n-1$ simply connected domains
$G_1,\dots,G_{n-1}$. Let $\zeta(u)$ be the function which is
defined on $G_k$ with respect to a fixed start point $\zeta_0$
with $\varphi(\zeta_0)=z_n$. Then the mappings
$\Phi_k:=\varphi|G_k=\varphi|G_k:G_k\to \{u\in
\overline{\C}\,:\,|u|>1\}$ are conformal.\par The boundaries of
the domains $G_j$ are pairwise disjoint. Each $\partial G_j$ is
mapped homeomorphically by $\varphi$ on the unit circle.\par It
holds $P(z,u):=\frac{Q(z,u)}{u}=(\frac{z}{u}-1)q(z)$. The
derivative zeros of $P$ with respect to $z$ are the same as those
of $Q$. For $u\to\infty$ the polynomials $P(z,u)$ tend locally
uniformly to $q(z)$. So, in this case, $\zeta(u)$ tends to
$\infty$ on one $G_k$, let us say on $G_1$. For $k=2,\dots,n-1$ it
follows that each $\zeta_j(u)\in G_k$ tends to some derivative
zero  $\xi_k$ of $q'$ if $u\to \infty$.
\subsection{$\zeta(u)$ on $G_1$}
From (\ref{fneu}) we see that $\zeta(u)$ has a pole in $\infty\in
G_1$. The equation
$$1-\frac{u}{\zeta(u)}=-\frac{q(\zeta(u))}{\zeta(u)q'(\zeta(u))}$$
gives that $\frac{u}{\zeta(u)}\to \frac{n+1}{n}$. It holds
$$f(u,\zeta(u))= \frac{z_n\, \zeta_0 \,q'(\zeta_0)}{q(\zeta_0)^2}
\cdot\frac{q(\zeta(u))^2}{u\,\zeta(u)\,q'(\zeta(u))}= \frac{z_n\,
\zeta_0 \,q'(\zeta_0)}{q(\zeta_0)^2}\frac{q(\zeta(u))}
{\zeta(u)\,q'(\zeta(u))}\frac{\zeta(u)}{u}
\frac{q(\zeta(u))}{\zeta(u)}.$$ All fractions stay to be finite
(and non zero) for $u\to\infty$, except of the last one, which has
a pole of order $n-2$ in $\infty$, and so $f$ has.
\subsection{$\zeta(u)$ on $G_k$ for $k>1$}
In this cases $\zeta(u)$ tends to some derivative zero $\xi_k$ of
$q'$.  From $$0=\zeta(u) q'(\zeta(u))- uq'(\zeta(u))+q(\zeta(u))$$
we conclude that $u q'(\zeta(u))\to q(\xi_k)$ if $u\to \infty$.
Now we see from (\ref{fneu}) that $f$ is holomorphic in
$\infty_k\in G_k$, and $c_k:=f(\infty_k)=\frac{z_n\, \zeta_0 \,
q'(\zeta_0)}{q(\zeta_0)^2}\cdot q(\xi_k)$. Thus $f$ is holomorphic
on $G_k$. Moreover $f$ does not vanish in $G_k$, because the zeros
of $q$ are all in $\overline{E}$. But on the boundary (as well as
on the boundary of $G_1$) there will be some zero, which comes
from the zero(s) of $p$ on the unit circle.
\section{Blowing up and pulling back}
Let $r>0$ and $p_r(z)=r^np(z/r)$. If we start the considerations
of the preceding section with $p_r$ instead of $p$ we have to
replace the zeros $z_2,\dots,z_n$ of $p$ by $rz_2,\dots,rz_n$ and
the derivative zeros $\zeta(u)$ by $r\zeta(u)$ as well as $q(z)$
by $r^{n-1}q(z/r)$. The variation is then
$$Q_r(z,u):=r^nQ(z/r,u)=z(z-ur)\cdot r^{n-1}q(z/r)=z(z-u
r)\prod_{j=2}^{n-1}(z-z_jr).$$ Note that the zeros of $Q_r(.,u)$
are the points $ru, rz_2,\dots,rz_n$, and it has the derivative
zeros $r\zeta(u)$, where $\zeta(u)$ denotes those of $Q(.,u)$.
\par  As already
mentioned we have $\rho(p_r)=\rho(p)$ for all $r>0$. Let $u_0$ be
some complex number of modulus $r$. If $r$ is large enough me may
provide that $$|f\big(ru_0(r),r\zeta(u_0(r))\big)|>|c_k|/2$$ if
$r\zeta(u_0(r))\in G_2,\dots,G_{n-1}$. If $r\zeta(u_0(r))\in G_1$
we may, because of the pole of $f$ in $\infty_1\in G_1$, assume
that $|f\big(ru_0(r),r\zeta(u_0(r))\big)|>1$. Now (\ref{fneu}) and
(\ref{Qf}) show $$\rho\big(Q_r(.,u),r\zeta(u)\big)=r\cdot f(ru,
r\zeta(u)).$$ So
$\rho\big(Q_r(.,u_0(r)),r\zeta(u_0(r))\big)>\rho(p,\zeta_0)$ for
all sufficiently large $r$ and all derivative zeros of this
polynomial.  If $\zeta_0$ has been taken above as an essential
derivative zero of $p$ this says that
$$\rho\big(Q_r(.,u_0(r)),\zeta(u_0(r))\big)>\rho(p)$$ for all
derivative zeros  $r\zeta(u_0(r))$ of $Q_r(.,u_0(r))$. This gives,
together with the remark above,
$\rho\big(Q_r(.,u_0(r))\big)>\rho(p_r)=\rho(p)$.\par The
polynomial $Q_r(.,u_0(r))$ has all its zeros in $|z|\le r$ and one
zero more on the boundary of this disk than $p$ have on the unit
circle (namely $u_0(r)$, in which $z_n$ has been changed). By
$p^*(z):=r^{-n}Q_r(z r,u_0(r))$ we pull all the zeros back into
the closed unit disk and so we found some polynomial, which has
one zero more on the unit circle as $p$ and which fulfills
$\rho(p^*)>\rho(p)$.\par We can repeat this argument until we
obtain a polynomial vanishing only on the unit circle and which
associated number is bigger than that of $p$. This finishes the
proof of theorem \ref{thm1}.
\section{Proof of Smale's conjecture}
It remains to compare $\rho(p)$ for polynomials
$p(z)=z\prod_{j=2}^n(z-z_j)$ with $|z_2|=\dots|z_n|=1$. For such
polynomials Smale's conjecture has already been proved by Tischler
\cite{tisch}. He also determined the maximal polynomials for this
subclass of $\Fn$ as $p(z)=a_1 z + a_n z^n$ with $a_1,a_n\in
\C\setminus\{0\}$, and we have the result that these are indeed
the only maximal polynomials in $\Fn$.

\bigskip

\begin{minipage}[t]{14cm}
Universit\"at Oldenburg\\ Fakult\"at V, Institut f\"ur
Mathematik\\Postfach 2503
\\ D-26111 Oldenburg\\ Bundesrepublik Deutschland\\ {\small
e-mail: schmieder@mathematik.uni-oldenburg.de}
\end {minipage}
\end{document}